\documentclass[12pt]{amsart}

\def\sw#1{{\sb{(#1)}}}
\def\sco#1{{\sb{<#1>}}} 
\def\su#1{{\sp{(#1)}}} 

\def\<{{\langle}}
\def\>{{\rangle}}

\def\eps{\epsilon}

\def\note#1{{}}

\def\note#1{}
\def\M{{\mathcal M}}

\def\cC{{\mathcal C}}

\def\ol{\overline}
\def\ot{\otimes}

\def\Label{\label}

\newtheorem{proposition}{Proposition}[section]

\newtheorem{corollary}[proposition]{Corollary}
\newtheorem{theorem}[proposition]{Theorem}

\theoremstyle{definition}
\newtheorem{definition}[proposition]{Definition}
\newtheorem{example}[proposition]{Example}

\theoremstyle{remark}
\newtheorem{remark}[proposition]{Remark}

\newcounter{c}

\newcommand{\etyk}[1]{\vspace{-7.4mm}$$\begin{equation}\Label{#1}
\addtocounter{c}{1}}
\renewcommand{\]}{\ifnum \value{c}=1 $$\else \end{equation}\fi}

\begin{document}

\title[Doi-Koppinen modules]{Doi-Koppinen modules for quantum
groupoids}
\dedicatory{To Max Kelly on the occasion of his 70th birthday.}
\author{Tomasz Brzezi\'nski}
\address{Department of Mathematics, University of Wales Swansea,
Singleton Park, Swansea SA2 8PP, U.K.}
\email{T.Brzezinski@swansea.ac.uk}
\urladdr{http//www-maths.swan.ac.uk/staff/tb}
\author{Stefaan Caenepeel}
\address{Faculty of Applied Sciences,
Free University of Brussels, VUB, Pleinlaan 2, B-1050 Brussels, Belgium}
\email{scaenepe@vub.ac.be}
\author{Gigel Militaru}
\address{Faculty of Mathematics, University of Bucharest, Str.
Academiei 14, RO-70109 Bucharest 1, Romania}
\email{gmilit@al.math.unibuc.ro}

\begin{abstract}
    A definition of a Doi-Koppinen datum over a noncommutative algebra
    is proposed. The idea is to
    replace a bialgebra in a standard Doi-Koppinen datum with a bialgebroid.
    The corresponding
    category of Doi-Koppinen modules over a noncommutative algebra
    is introduced. A weak Doi-Koppinen
    datum and module of \cite{Boh} are shown to be
    examples of a Doi-Koppinen datum and module over an algebra.  A coring
associated to a Doi-Koppinen
    datum over an algebra is constructed and various properties of
    induction and forgetful functors for Doi-Koppinen modules over an
    algebra are deduced from the properties of corresponding functors
    in the category of comodules of a coring.
\end{abstract}
\maketitle

\section{Introduction}
Doi-Koppinen modules introduced in \cite{Doi:uni} \cite{Kop:var} as a
generalisation of (co)modules or Hopf modules studied in Hopf algebra
theory can be viewed
as a representation of a triple comprising an algebra, a coalgebra
and a bialgebra which satisfy certain compatibility conditions. Recently
these have been generalised to the case in which a bialgebra is
replaced by a weak Hopf algebra \cite{Boh}. It is known
\cite{EtiNik:dyn} that weak Hopf algebras are an example of a
generalisation of a bialgebra known as an $R$-bialgebroid \cite{Lu}
or $\times_{R}$-bialgebra \cite{Tak:gro} (and leading to the notion of
a Hopf algebroid or a quantum groupoid), introduced in the context of
Poisson geometry, algebraic topology and
classification of algebras. It seems therefore natural to
ask whether a definition of a Doi-Koppinen datum in which a
bialgebra is replaced by an
$R$-bialgebroid is possible. In this paper we propose such a
definition and by this means introduce the notion of a
Doi-Koppinen module over a noncommutative algebra $R$.
We show that Doi-Koppinen
modules for a weak Hopf algebra are a special case thus providing a
new, more general point of view on weak Doi-Koppinen data and modules.

On the other hand it has been realised in \cite{Brz:mod} that a
natural point of view on Doi-Koppinen data is provided by entwining
structures introduced in \cite{BrzMaj:coa}.
The same point of view was adopted in \cite{CG}, where
weak entwining structures were introduced in order to describe
Doi-Koppinen data for a weak Hopf algebra. Later on it has been shown in
\cite{Brz:cor} that both entwined modules and weak entwined modules are
simply comodules of certain corings. Thus various properties of entwined
modules such as  Frobenius and separability properties discussed first
in the case of Doi-Koppinen modules in \cite{CIMZ} and \cite{CMZ},
 can be derived from the properties of comodules over a coring. In the
present paper we show that a Doi-Koppinen datum for an $R$-bialgebroid leads
to a certain coring whose comodules are precisely the Doi-Koppinen
modules over $R$.

\section{Preliminaries}
\subsection{Notation}\label{subs.notation}
We use the following conventions. For an object $V$ in a category,
the identity morphism $V\to V$ is denoted by $V$. All rings in this
paper have 1, a ring map is assumed to respect 1, and all
modules over a ring are assumed to be unital.
For a ring $R$,  $\M_R$ (resp.\ ${}_R\M$, ${}_{R}\M_{R}$) denotes
the category of right $R$-modules (resp.\ left $R$-modules,
$R$-bimodules).
The action of $R$ is denoted by a dot between elements.

Throughout the paper $k$ denotes a commutative ring with unit.
We  assume that all the algebras are over $k$ and  unital, and
coalgebras are over $k$ and  counital. Unadorned tensor product
is over $k$.
For a
$k$-coalgebra $C$ we use $\Delta_{C}$ to denote the coproduct and $\eps_{C}$ to
denote the counit (we skip subscripts if no confusion is possible).
Notation for comodules is similar to that for modules
but with subscripts replaced by superscripts, i.e. ${}^C\M$ is the
category of left $C$-comodules etc. We
use the Sweedler notation for coproducts and coactions, i.e. $\Delta(c)
= c\sw 1\otimes c\sw 2$ for a coproduct, and
$\rho(m) = m_{<-1>} \otimes m_{<0>}$ for a left coaction (summation
understood).

Let $R$ be a $k$-algebra. Recall from \cite{Swe:pre} that an
{\em $R$-coring} is a coalgebra in the monoidal
category of $R$-bimodules $({}_R\M_R, \otimes_R , R)$, i.e.,
it is an $(R,R)$-bimodule
$\cC$ together with $(R,R)$-bimodule maps $\Delta_\cC:\cC\to
\cC\otimes_R\cC$ called a coproduct and $\eps_\cC:\cC\to R$ called a
counit, such that
$$
(\Delta_\cC\otimes_R\cC)\circ\Delta_\cC = (\cC\otimes_R\Delta_\cC)\circ
\Delta_\cC, \quad (\eps_\cC\otimes_R\cC)\circ\Delta_\cC =
(\cC\otimes_R\eps_\cC)\circ \Delta_\cC = \cC.
$$
We use the Sweedler notation for the coproduct $\Delta_{\cC}$ too.
A  left $R$-module $M$
together with a left
$R$-module map ${}^M\rho:M\to \cC \otimes_R M$
such that
$$
(\cC \otimes_R \;{}^M\rho)\circ \rho = (\Delta_\cC \otimes_R M)\circ
{}^M\rho, \quad
(\eps_\cC \otimes_R M)\circ {}^M\rho = M
$$
is called a {\em left comodule of the coring $\cC$} or, simply, a {\em left
$\cC$-comodule}, and ${}^M\rho$ is called a {\em left coaction}.
A map between left $\cC$-comodules is a left $R$-module map
$f:M\to N$ such that
${}^N\rho\circ f = (\cC \otimes_R f)\circ {}^M\rho$.
The category of
left $\cC$-comodules is denoted by ${}_{R}^{\cC}\M$.

\subsection{$R$-rings and bialgebroids}
Let $R$ be a $k$-algebra. Recall from \cite{Swe:alg}, \cite{Tak:gro}
that an {\em $R$-ring} is a pair
$(U, i)$, where $U$ is a $k$-algebra  and $i: R\to U$ is an algebra map.
If $(U, i)$ is an $R$-ring then $U$ is an $(R,R)$-bimodule  with the
structure provided by the map $i$,
$r\cdot u \cdot r':= i(r)ui(r')$.

Let $R$ be an algebra and $\bar{R}=R^{op}$ the opposite algebra, and
let $R^{e}=R\otimes \bar{R}$ be the enveloping algebra of $R$.
In case $(H,i)$ is an $R^{e}$-ring, the map $i$ is necessarily of
the form $i=m_H\circ(s_{H}\otimes t_{H})$, where
 $s_H: R \to H$, $t_H: \bar{R}\to H$ are algebra maps such that
$s_H(a)t_H(\bar{b})=t_H(\bar{b}) s_H(a)$,
for all $a\in R$, $\bar{b}\in \bar{R}$, and $m_H$ is the product in $H$.
In this case   $s_H$ is called the {\em source}
map and $t_H$ the {\em target} map.
$(H,i=m_H\circ(s_H\otimes t_H))$ is denoted by $(H,s_H,t_H)$.

Let $(H, s_H, t_H)$ be an $R^{e}$-ring and $(A, s_A)$ an $R$-ring.
We view $H$ as an  $R$-bimodule,
via $r\cdot h \cdot r' = s_H(r)t_H(r')h$,
and $A$ as an $R$-bimodule via $s_A$, and define \cite{Tak:gro}
$$
H\times_R A = \{\sum_{i} h^{i}\otimes_{R} a^{i}\in H\otimes_R A \; |\;
\forall r\in R ,\;
\sum_{i} h^i t_H(r) \otimes_R a^i = \sum_{i} h^i\otimes_R a^i s_A(r)\}.
$$
$H\times_R A$ is an $R$-ring
with product
$$
(\sum_{i} h^i\otimes_R a^i) (\sum_{j} \tilde{h}^j\otimes_R \tilde{a}^j)=
\sum_{i,j} h^i \tilde{h}^j\otimes_R a^i\tilde{a}^j,
$$
the unit $1_H\otimes_R 1_A$ and the algebra map
$R\to H\times_R A$, $a\mapsto s_H (a)\otimes_R 1_A$
(cf. \cite{Tak:gro}). Taking $(A, s_A) = (H, s_H)$ we can define
$H\times_R H$ which is not only an $R$-ring but also an $R^e$-ring via
$R\otimes \bar{R} \to H\times_R H$,
$a \otimes_R \bar{b}\mapsto s_H(a)\otimes_R t_H(\bar{b})$.

\begin{definition}\label{def.Lu}
Let $(H, s_H, t_H)$ be an $R^{e}$-ring. We say that $(H, s_H, t_H,
\Delta, \eps)$ is
an {\em $R$-bialgebroid} iff $(H,\Delta, \eps)$ is an $R$-coring such
that
${\rm Im}(\Delta) \subseteq H\times_R H$ and the corestriction
of the coproduct $\Delta :H \to H\times_R H$ is an algebra map,
$\eps (1_H)=1_R$, and for all $g$, $h\in H$
\begin{equation}\label{25}
\eps(gh) =\eps \Bigl( g s_H( \eps (h) ) \Bigl) =
\eps \Bigl( g t_H( \eps (h) ) \Bigl).
\end{equation}
\end{definition}
It is shown in \cite{BM} that this is equivalent both to the definition
of a bialgebroid in \cite{Lu} and that of $\times_R$-bialgebra in
\cite{Tak:gro}.

Szlach\'anyi \cite{Szl} has reformulated the definition of bialgebroid in
terms of monoidal categories and monoidal functors: if $H$ is an
$R^e$-ring, then we have the restriction of scalars functor
$F:\ {}_H{\mathcal M}\to{}_R{\mathcal M}_R$. $H$ is an $R$-bialgebroid
if and only if there exists a monoidal structure on ${}_H{\mathcal M}$
such that $F$ is a strict monoidal functor. If $(H, s_H, t_H,
\Delta, \eps)$ is as in Definition \ref{def.Lu}, then the corresponding
monoidal structure
on ${}_H{\mathcal M}$ is given by
$$h\triangleright (m\otimes_R n)=h_{(1)}m\otimes_R h_{(2)}n~~;~~
h\triangleright a=\varepsilon(hs(a))=\varepsilon(ht(a))$$
for all $m\in M\in {}_H{\mathcal M}$, $n\in N\in {}_H{\mathcal M}$, $a\in R$.

Basic examples of $R$-bialgebroids are provided by $R^e$ and
 ${\rm End}(R)$, in the
case
 $R$ is finitely
generated projective over $k$ (see \cite{Lu}, \cite{Tak:gro}). In
particular
any matrix algebra $M_n (k)$ has a structure of an
$R$-bialgebroid
with an antipode over any $n$-dimensional algebra $R$. We believe,
this gives a nice motivation for studying bialgebroids from
an algebraic point of view.

\subsection{Doi-Koppinen datum over a weak Hopf algebra}
A {\em weak bialgebra} is an algebra and a coalgebra $H$ with
multiplicative (but non-unital) coproduct such that for
all $x,y,z\in H$,
$\eps(xyz) = \eps(x y\sw 1)\eps(y\sw 2 z) =
\eps(x y\sw 2)\eps(y\sw 1 z)$,
and
\begin{equation}
(\Delta\otimes H)\circ \Delta(1) = (\Delta(1)\otimes 1)(1\otimes
\Delta(1)) = (1\otimes \Delta(1))(\Delta(1)\otimes 1).
\label{weak.unit}
\end{equation}
A  {\em weak Hopf algebra} is  a weak bialgebra $H$ with an antipode,
i.e., a
linear map $S:A\to A$ such that for all $h\in H$,
$h\sw 1 S(h\sw 2) = \eps(1\sw 1h)1\sw 2$, $S(h\sw 1)h\sw 2 = 1\sw
1\eps(h 1\sw 2)$, and  $S(h\sw 1) h\sw 2 S(h\sw 3) = S(h)$.
Weak Hopf algebras have been introduced
in \cite{BohSzl:coa} \cite{Nil:axi}
and studied  in connection to
integrable models and classification of subfactors of von Neumann
algebras. Given a weak Hopf algebra $H$ with bijective antipode,
define the maps,
$$
\Pi^L, \Pi^R :H\to H, \quad
\Pi^L(g)=\eps (1_{(1)}g)1_{(2)}, \quad
\Pi^R(g)=\eps (g 1_{(2)})1_{(1)}.
$$
 Then \cite{Boh:wea}
$R :={\rm Im} (\Pi^L)$ is a subalgebra of $H$, separable and
 Frobenius as a $k$-algebra with the separability
idempotent $e = S(1_{(1)}) \otimes 1_{(2)}\in R\otimes R$ and the
Frobenius pair $(e, \varphi)$, where  $\varphi := \eps{|_{R}}$.
The fact that $e$ is a separability idempotent means explicitly
\begin{equation}\label{sepide}
\forall g\in H, \quad \Pi^L (g) S(1_{(1)}) \otimes 1_{(2)} =
S(1_{(1)}) \otimes 1_{(2)} \Pi^L (g).
\end{equation}
Numerous useful formulae for a weak Hopf
algebra were proven in  \cite{Boh:wea}.
Although some of them, such as
(\ref{sepide}), were obtained using duality arguments valid only in the
finite dimensional case, one can also prove them directly using the
axioms  of a weak Hopf algebra. The proofs are not always obvious,
but quite straightforward once one becomes familiar with these axioms.

Finally, recall from \cite{Boh} the following

\begin{definition}\label{def.weakDK}
A left-left {\em weak Doi-Koppinen datum} is
a triple $(H, A, C)$, where $H$  is a weak Hopf algebra and

(1) $(A, {}^A\rho)$ is a left {\em weak $H$-comodule algebra},
i.e., $A$ is
an algebra and a left $H$-comodule such that
${}^A\rho (a){}^A\rho (b)={}^A\rho (ab)$, and
$(H\otimes {}^A\rho)\circ{}^A\rho (1)= \sum 1_{(1)}\otimes
1_{<-1>}1_{(2)}\otimes 1_{<0>}$,
for all $a,b\in A$;

(2) $C$ is a left {\em weak $H$-module coalgebra}, i.e.,
$C$ is a coalgebra with counit $\eps_C$
and a left $H$-module such that
$\Delta_{C}(h\cdot c) =
\sum h_{(1)} \cdot c_{(1)}\otimes h_{(2)} \cdot c_{(2)}$, and
$\eps_C (hg\cdot c)= \eps_H (hg_{(2)})\eps_C (g_{(1)} \cdot c)$
for all $c\in C$ and $h,g\in H$.
\end{definition}

A  (left-left) {\em weak Doi-Koppinen module} associated to a weak
Doi-Koppinen datum $(H,A,C)$ is a triple
$(M, \cdot, {}^M\rho)$,
where $(M, \cdot)$ is a left $A$-module, $(M, {}^M\rho)$ is a left
$C$-comodule, and
$$
{}^M\rho (a\cdot  m) = a_{<-1>} \cdot m_{<-1>} \otimes a_{<0>}
\cdot m_{<0>}.
$$
Note that here ${}^A\rho(a) = a\sco{-1}\otimes a\sco 0 \in H\otimes C$
and ${}^M\rho(m) = m\sco{-1}\otimes m\sco 0 \in C\otimes M$.  The
category of (left-left) weak Doi-Koppinen modules is denoted by
${}^C_A \M(H)$.

Morphisms between left weak $H$-comodule algebras (resp. left weak
$H$-module coalgebras) are defined in the obvious way: they are
$k$-linear maps that are $H$-colinear (resp. $H$-linear) algebra (resp.
coalgebra) maps. Thus we can consider the categories of left
weak $H$-comodule algebras, left weak $H$-module coalgebras and
left-left weak Doi-Koppinen data over $H$. The latter is denoted by
$\mathcal W\mathcal D\mathcal K(H)$.

\section{Doi-Koppinen modules over algebras}

In this section we define the notion of a Doi-Koppinen datum over
a noncommutative algebra and we relate it to a weak Doi-Koppinen datum.
Our definition is in part
motivated by  the following important observation
\cite[Proposition~2.3.1]{EtiNik:dyn}.

\begin{proposition} \label{webial}
Let $H$ be a  weak Hopf algebra with coproduct $\Delta$, counit $\eps$,
and bijective antipode $S$,
and let $R ={\rm Im} (\Pi^L)$. Then
$H$ is an $R$-bialgebroid
with the source and target $s_H$, $t_H :R \to H$ given by
$$
s_H (\Pi^L (g) ) = \Pi^L (g), \quad
t_H (\Pi^L (g) ) = S^{-1} ( \Pi^L (g) ) = \eps (1_{(2)} g) 1_{(1)},
$$
and the comultiplication $\tilde{\Delta} : H\to H\otimes_R H$ and the
counit $\tilde{\eps} :H \to R$ given by
$$
\tilde{\Delta} (h) = (can \circ \Delta)(h) = h_{(1)} \otimes_R h_{(2)}, \quad
\tilde{\eps} (h) = \Pi^L (h)
$$
for all $h\in H$, where
$can : H\otimes H \to H\otimes _R H$ is the canonical projection.
\end{proposition}

\begin{proof}
For the details we refer to \cite{EtiNik:dyn}, we only remark
that  ${\rm Im} (\tilde{\Delta}) \subseteq H\times_R H$ can be
established from the separability of $R$ as follows.
Apply $S^{-1} \otimes H$  to (\ref{sepide}),
for an arbitrary $h\in H$ write $h=h1_H$, and use that $\Delta$ is
multiplicative
to obtain
\begin{equation}\label{2.31a}
h_{(1)} \otimes h_{(2)} \Pi^L (g) = h_{(1)} 1_{(1)} \otimes
h_{(2)} 1_{(2)} \Pi^L (g) = h_{(1)} S^{-1} (\Pi^L (g) ) \otimes h_{(2)}.
\end{equation}

We also note that
in  \cite{EtiNik:dyn} the conditions (\ref{25}) are not
required for an $R$-bialgebroid. However it can be easily seen that
$\tilde{\eps}$ as defined above satisfies equations (\ref{25}).
\end{proof}

There is also a partial converse to Proposition~\ref{webial} (cf.\
\cite[Proposition 1.6]{Szl}).
Suppose $H$ is an $R$-bialgebroid, where $R$ is a Frobenius
and separable  $k$-algebra. Let $e=e\su 1 \otimes e\su 2$ (summation
understood) be a
separability idempotent and let
$\varphi :R \to k$ be the Frobenius map such that $(e,\varphi)$ is a
Frobenius pair. Then $H$ is a weak bialgebra
with the coproduct $\tilde{\Delta}: H\to H\otimes H$,
$h\mapsto h_{(1)} \cdot e\su 1 \otimes e\su 2 \cdot h_{(2)}$, and the counit
$\tilde{\eps} = \varphi \circ \eps  :H\to k$, where $\eps:H\to R$ is the
counit of the $R$-coring $H$.

\begin{definition}\label{altadef}
Let $(H, s_H, t_H)$ be an $R$-bialgebroid. Then a
{\em left $H$-module coalgebra} is a coalgebra $C$ in the monoidal
category $(_H\M, \otimes_R, R)$ of left $H$-modules.
\end{definition}

Recall from \cite{Szl} that ${}_H\M$ has a monoidal structure defined
as follows.
For all $M, N\in {}_H\M$, $M\otimes_R N\in {}_H\M$ via
$h\cdot (m\otimes_R n) = h_{(1)}\cdot m\otimes_R h_{(2)}\cdot n$.
$R$ is the unit object, when viewed in ${}_H\M$ via the action
\begin{equation}\label{Stef2}
h\triangleright a = \eps (h s_H(a))= \eps (h t_H(a)).
\end{equation}

Thus, $C$ is a left $H$-module coalgebra
if and only if $(C, \cdot)$ is a left $H$-module and $(C, \Delta_C,
\eps_C)$ is an
$R$-coring, where $C$ is viewed as an
 $R$-bimodule via
$r\cdot c \cdot r'= s_H(r) t_H(r') \cdot c$, such that
$\Delta_C$, $\eps_C$ are left
$H$-modules maps, i.e., for all $h\in H$ and $c\in C$,
\begin{equation}\label{Stef3}
\Delta_C (h\cdot c) = h_{(1)} \cdot c_{(1)} \otimes_R
h_{(2)} \cdot c_{(2)}, \quad
\eps_C (h\cdot c) =  h\triangleright \eps_C (c) =
\eps_H (h s_H (\eps_C (c)) ).
\end{equation}
A morphism between two $H$-module coalgebras is an
 $H$-linear map of $R$-corings. We can then
consider the category of $H$-module coalgebras.

\begin{example}
(1) $(H, \Delta_H, \eps_H)$ is a left $H$-module coalgebra
with the left $H$-action provided by the left multiplication.

(2) View $R$ as an $R$-coring in the trivial way, i.e., both $\Delta_R$ and
$\eps_R$ are identity maps. Then
$(R, \triangleright )$ is a left $H$-module coalgebra.
Indeed, note  that
for all $r\in R$ and $h\in H$ we have that
$s_H(h_{(1)}\triangleright r)h_{(2)}=hs_H(r)$ \cite{BM},
and then apply the left $R$-module map $\eps_H$ to  obtain that
$\eps_H (h s_H (r)) = \eps_H (h_{(1)} s_H(r) )\eps_H (h_{(2)})$. This
is equivalent to the fact that $\Delta_R$ is a left $H$-linear
map.

(3) $C=R^e$ is an $R$-coring with the coproduct $\Delta_{R^e}
(r \otimes \bar{r}) = r \otimes 1_{\bar{R}} \otimes_R 1_{R} \otimes
\bar{r}$ and the counit $\eps_{R^e} (r \otimes \bar{r}) = r\bar{r}$, and
it can be made into a left $H$-module coalgebra by the $H$-action
$h\cdot (r \otimes \bar{r}) = \eps_H (h s_H(r) t_H(\bar{r}) ).$
\end{example}

\begin{definition}\label{lhca}
Let $(H,s_H,t_{H})$ be an $R$-bialgebroid. A {\em left $H$-comodule
algebra} is a triple
$(A, s_A, {}^A\rho)$ where

(1) $(A, s_A)$ is an $R$-ring.

(2) $(A, {}^A\rho)$ is a left comodule of the $R$-coring $H$.

(3) ${\rm Im} ({}^A\rho) \subseteq H\times_R A$ and its corestriction
${}^A\rho : A \to H\times_R A$ is an algebra map.
\end{definition}

A morphism between two $H$-comodule algebras is a left $H$-colinear map
that is also a morphism of $R$-rings. A morphism of $R$-rings is
defined in the obvious way. Thus we can consider the category of
$H$-comodule algebras.

\begin{example}
(1) $(H, s_H, \Delta)$ is a left $H$-comodule algebra.

(2) $(R, s_R=R, {}^R\rho)$, where  $
{}^R\rho :R \to H\otimes_R R$,
$r\mapsto s_H(r) \otimes_R 1_R$ is a (trivial)
 left $H$-comodule algebra

(3) $A=R^e$ is a left $H$-comodule algebra via
$s_{R^e} = R\otimes 1_R$, and
${}^{R^e}\rho (r \otimes \bar{r}) = s_H (r) \otimes_R 1_R \otimes
\bar{r}$.

(4) Example~(3) can be generalised as follows.
For an $H$-comodule algebra $A$ and an algebra $B$,
$A\otimes B$ is a left $H$-comodule algebra with the structures arising
from those of $A$.

(5) An interesting nontrivial example of a comodule algebra of a
bialgebroid is constructed in \cite[Theorem~6.3, Lemma~6.7]{Sch:bia}.
Let $H$ be a Hopf algebra and let $A/B$ be an $H$-Galois extension,
i.e., $A$ be a right $H$-comodule algebra with a right coaction
$\rho^A:A\to A\otimes H$, $a\mapsto a\sco 0\otimes a\sco 1$, $B = A^{co
H}$, and the canonical map $\chi:A\otimes_B A\to A\otimes H$, $a\otimes_B
a'\mapsto aa'\sco 0\otimes a'\sco 1$ be bijective. Suppose $H$ is
$k$-flat and $A$ is a faithfully flat left $B$-module. View $A^e$ as a
right $H$-comodule via $a\otimes \bar{a} \mapsto a\sco 0\otimes
\bar{a}\sco 0\otimes a\sco 1\bar{a}\sco 1$. Then the $k$-module of
coinvariants $G=(A^e)^{co H}$ is a subalgebra of $A^e$ and a
$B$-bialgebroid via $s_G: b\mapsto b\otimes 1$, $t_G: b\mapsto 1\otimes
b$, $\Delta_G: \sum_i a^i\otimes \bar{a}^i \mapsto  \sum_ia^i\sco 0\otimes
\chi^{-1}(1\otimes a^i\sco 1)\otimes \bar{a}^i$, and $\eps_G : \sum_i
a^i\otimes \bar{a}^i \mapsto \sum_ia^i\bar{a}^i$. Furthermore $A$ is a left
comodule algebra of the $B$-bialgebroid $G$ with a left coaction
${}^A\rho: A\to G\otimes_B A$, $a\mapsto a\sco 0\otimes
\chi^{-1}(1\otimes a\sco 1)$.
\end{example}

Note that the definition of a left $H$-comodule algebra is not dual to
that of a left $H$-module coalgebra. The reason is that, although the
category of left $H$-modules is monoidal, the category ${}^H_R \M$
of left comodules of an $R$-coring $H$ is not.
Thus there is no
way of defining a left $H$-module algebra as an algebra in the
category ${}^H_R\M$.
However one can consider more restrictive definition of a left
$H$-comodule $M$ (cf.\ \cite[Definition~5.5]{Sch:bia}) by requiring it to
be an $R$-bimodule with
an $R$-bimodule coaction
${}^M\rho:M\to H \otimes_R M$ such that ${\rm Im}({}^M\rho)
\subseteq H\times_R M$,
where
$$
H\times_R M =
\{ \sum_{i} h^{i}\otimes_{R} m^{i}\in H\otimes_R M \; |\; \forall r\in R ,\;
\sum_{i} h^i t_H(r) \otimes_R m^i = \sum_{i} h^i\otimes_R m^i \cdot r\}.
$$
The subcategory ${}^H_R \M_R \subseteq {}^H_R \M$ of all such comodules
is monoidal.  For all $M$, $N\in {}^H_R \M_R$,
$M\otimes_R N \in {}^H_R \M_R$ via
${}^{M\otimes_RN}\rho (m\otimes_R n) = m_{<-1>} n_{<-1>} \otimes_R m_{<0>}
\otimes_R
n_{<0>}$.
Note that the right hand side is well defined because
${\rm Im}({}^{M}\rho) \subseteq H\times_R M$.
$R$ is the unit in ${}^H_R \M_R$  with the trivial coaction
${}^R\rho (r) = s_H (r) \otimes_R 1_R$. Furthermore
the forgetful functor $F: {}^H_R \M_R \to {}_R\M_R$ is strict monoidal.
Now, one could define a left $H$-comodule algebra as an algebra
in the monoidal category
$({}^H_R \M_R, \otimes_R, R)$. It appears, however, that this definition
is too restrictive to cover the case of a weak Doi-Koppinen datum.

\begin{definition}\label{dhdat}
Let $(H, s_H, t_H, \Delta, \eps)$ be an $R$-bialgebroid. Then
$(H, A, C)$ is called a {\em (left-left) Doi-Koppinen datum over
(an algebra) $R$} if
$A$ is a left $H$-comodule algebra and $C$ is a left
$H$-module coalgebra. Such a datum is denoted by $(H,A,C)_{R}$. The
category of Doi-Koppinen data with $H$ over $R$ is denoted by $\mathcal
D\mathcal K_R(H)$.

A {\em (left-left) Doi-Koppinen module} over $R$ (associated to
$(H,A,C)_{R}$) is a triple
$(M, \cdot, {}^M\rho)$, where $(M, \cdot)$ is a left $A$-module
(hence $M$ is a left $R$-module via $s_A$),
$(M, {}^M\rho)$ is a left comodule of the $R$-coring $C$,
and for all $a\in A$ and $m\in M$,
\begin{equation}\label{comp}
{}^M\rho (a \cdot m) = a_{<-1>} \cdot m_{<-1>} \otimes_R a_{<0>} \cdot m_{<0>}.
\end{equation}
\end{definition}

The category of Doi-Koppinen modules associated to $(H, A, C)_{R}$
is denoted by  ${}^C_A\M (H;R)$.
Note that the right hand side of (\ref{comp})
is well defined since
${\rm Im} ({}^A\rho) \subseteq H\times_R A$.

\begin{example}
There are various examples of special
  cases of the category ${}^C_A\M (H;R)$ obtained by setting
$A = H, R, R^e$ and $C = H, R, R^e$. In particular, the category
of left $H$-modules, the category of left $H$-comodules or
the category of (generalised) relative Hopf modules
${}^H_A\M (H;R)$ and its
dual
${}^C_H\M (H;R)$ are all special cases of the category ${}^C_A\M (H;R)$.
\end{example}

The main aim of this section is to show that a weak
Doi-Koppinen datum in Definition~\ref{def.weakDK}
 is a
special case of a Doi-Koppinen datum over a noncommutative algebra. This
provides one with a new point of view on weak Doi-Koppinen modules.
In the proof of the next two Propositions, we will make use of the
following remark.

\begin{remark}\label{separableremark}
If $S$ is a
separable $k$-algebra with an idempotent $e= e\su 1 \otimes e\su 2$
(summation understood), and
$M$ and $N$ are $S$-bimodules then the canonical
projection $M\otimes N \to M\otimes_S N$ has a section
$$
\sigma: M\otimes_S N \to M\otimes N, \quad
\sigma(m\otimes_S n) = m\cdot e\su 1 \otimes e\su 2 \cdot n.
$$
If $H$ is a weak Hopf algebra, then $H$ is a bialgebroid over a separable
(and Frobenius) algebra $R$ with idempotent
$e = S(1_{(1)}) \otimes 1_{(2)}$ and the Frobenius map
$\varphi :R\to k$,  $\varphi := \eps{|_{R}}$, the restriction of a weak
counit of $H$ to $R$.
\end{remark}

\begin{proposition}\label{comodulealgebra}
Let $H$ be a weak Hopf algebra viewed as an $R$-bialgebroid
as in Proposition \ref{webial}. Then the category of left weak
$H$-comodule algebras (in the sense of
Definition~\ref{def.weakDK}) and of
left comodule algebras over the bialgebroid $H$ (in the sense of
Definition~\ref{lhca}) are isomorphic to each other.
\end{proposition}

\begin{proof}
1) Let $A$ be a left $H$-comodule algebra. We will show that $A$ is
then a left comodule algebra over the $R$-bialgebroid $H$.
Note that
$$s_A:\ R \to A, \quad s_A (\Pi^L (h)) =\eps (1_{<-1>} h ) 1_{<0>},$$
where $\eps:H\to k$ is a weak counit,
is a well-defined algebra map. Indeed, suppose $r=\Pi^L (h)=0$. This
means that $\eps (1_{(1)} h ) 1_{(2)} =0$, and therefore
$$
s_A (r) = \eps (1_{<-2>} h) \eps (1_{<-1>} ) 1_{<0>}
=
\eps (1_{(1)} h ) \eps (1_{(2)} 1_{<-1>} ) 1_{<0>} =0.
$$
Note that the second equality was obtained by using the following
observation made in \cite[Proposition~4.11]{CG}. The unit property
of a coaction of a weak comodule algebra in
Definition~\ref{def.weakDK}(1) is equivalent to
\begin{equation}\label{64}
1\sco{-2}\otimes 1\sco{-1}\otimes 1\sco 0 =
1_{(1)}\otimes 1_{(2)}1_{<-1>}\otimes 1_{<0>}.
\end{equation}
This proves that $s_A$ is well-defined. To prove that $s_A$ is an
algebra map we require the following two equalities (cf.\ \cite[Eq.\
(2.9b)]{Boh:wea} and \cite[Eq.\ (2.1b)]{Boh} respectively). For all
$g,h\in H$ and $a\in A$,
\begin{equation}\label{foreign}
\Pi^R(g)h = h\sw 1\eps(gh\sw2), \qquad {\Pi}^R(a_{<-1>})\otimes
a_{<0>}=1_{<-1>}\otimes a1_{<0>}.
\end{equation}
Now for
all
$r=\Pi^L (h)$, $s=\Pi^L (g)$ we have
\begin{eqnarray*}
s_A (rs) &=& s_A (\Pi^L (h) \Pi^L (g) )\\
&=& s_A \Bigl(\Pi^L (\Pi^L (h) g  )  \Bigl)
= \eps \Bigl( 1_{<-1>} \eps (1_{(1)} h ) 1_{(2)} g  \Bigl) 1_{<0>} \\
&=& \eps (1_{<-2>} h) \eps (1_{<-1>} g) 1_{<0>}
= \eps (1_{<-2>} h) \eps(g_{(1)}) \eps (1_{<-1>} g_{(2)}) 1_{<0>} \\
~~\mbox{(by Eqs~(\ref{foreign}))}~~
&=& \eps (1_{<-2>} h) \eps (\Pi^R (1_{<-1>}) g) 1_{<0>}
=\eps (1_{<-1>} h) \eps (1_{<-1'>} g) 1_{<0>}1_{<0'>}\\
&=& s_A (r) s_A (s),
\end{eqnarray*}
where  $1\sco{-1'}\otimes 1\sco{0'}$ denotes another copy of
$1\sco{-1}\otimes 1\sco 0$, and \cite[Lemma~2.5]{Boh:wea} has been used to
derive  the second equality and
the unit property in Definition~\ref{def.weakDK}~(1) to obtain the fourth
one. This proves
that $s_A$ is an algebra map and hence $(A, s_A)$ is an $R$-ring.

Next, using the canonical
projection $can : H\otimes A \to H\otimes_R A$ define a map
$\tilde{\rho} = can \circ {}^A\rho$, $\tilde{\rho}: A\to H\otimes_R
A$, where ${}^A\rho: A\to H\otimes A$ is the left coaction of a weak
Hopf algebra.
Explicitly,
$\tilde{\rho} (a) = a_{<-1>} \otimes_R a_{<0>}$. The map
$\tilde{\rho}$ is
left $R$-linear since using equation (\ref{64}) and the fact that $A$ is
a comodule algebra we have
for all $r=\Pi^L(h)\in R$ and $a\in A$
\begin{eqnarray*}
\tilde{\rho}(r\cdot a) & = &\tilde{\rho}(\eps(1\sco{-1}h)1\sco 0a) =
\eps(1\sco{-2}h)1\sco{-1}a\sco{-1}\otimes_R1\sco 0 a\sco 0 \\
&=& \eps(1\sw 1h)1\sw 2
1\sco{-1}a\sco{-1}\otimes_R 1\sco 0 a\sco 0 = \Pi^L(h)a\sco{-1}\otimes_R
a\sco 0.
\end{eqnarray*}
Then it is clear that
$(A, \tilde{\rho}) \in {}^H_R\M$.  We
prove now that ${\rm Im}(\tilde{\rho} )\subseteq H\times_R A$.
Equation (\ref{2.31a}) implies that for all $g\in H$ and $a\in A$
$$
a_{<-2>} \otimes a_{<-1>} \Pi^L (g) \otimes a_{<0>} =
a_{<-2>} S^{-1} (\Pi^L (g)) \otimes a_{<-1>} \otimes a_{<0>}.
$$
Apply $H\otimes \eps \otimes A$ to the last equality to obtain
$$
a_{<-2>} \otimes \eps (a_{<-1>} \Pi^L (g) ) a_{<0>} =
a_{<-1>} S^{-1} (\Pi^L (g)) \otimes a_{<0>}
$$
Using equations (2.2a),(2.2b) in \cite{Boh:wea} which, put together,
 state that for
all $g,h\in H$, $\eps(g\Pi^L(h)) =\eps(\Pi^R(g)h)$
we compute
\begin{eqnarray*}
a_{<-2>} \otimes \eps (a_{<-1>} \Pi^L (g) ) a_{<0>} &=&
a_{<-2>} \otimes \eps (\Pi^R (a_{<-1>}) g ) a_{<0>}\\
~~\mbox{ (by Eq.\ (\ref{foreign}))}~~
&=& a_{<-1>} \otimes \eps (1_{<-1>} g) a_{<0>} 1_{<0>} \\
&=& a_{<-1>} \otimes a_{<0>} s_A (\Pi^L (g)).
\end{eqnarray*}
Hence we have proved that
\begin{equation}\label{eq.sep}
    a_{<-1>} S^{-1} (\Pi^L (g)) \otimes a_{<0>} =
a_{<-1>} \otimes a_{<0>} s_A (\Pi^L (g)).
\end{equation}
In particular,
${\rm Im}(\tilde{\rho}) \subseteq H\times_R A$ as required.
It remains to be proven that $\tilde{\rho} (1_A) = 1_H \otimes_R 1_A$.
Using the unit property of a comodule algebra of a weak Hopf algebra in
Definition~\ref{def.weakDK}(1)
 we compute
$$
\rho (1_A) =
1_{<-2>} \otimes \eps (1_{<-1>}) 1_{<0>}
= 1_{(1)} \otimes \eps (1_{<-1>} 1_{(2)} ) 1_{<0>} = 1_{(1)} \otimes
s_A(1_{(2)} ).
$$
Since $A$ is a left $R$-module via $s_A$ we obtain
$$
1_{<-1>} \otimes_R 1_{<0>} = 1_{(1)} \otimes_R s_A(1_{(2)} ) =
1_{(1)}\cdot 1_{(2)} \otimes_R 1_A = S^{-1} (1_{(2)}) 1_{(1)}\otimes_R 1_A =
1_H \otimes_R 1_A.
$$
This completes the proof that $(A, s_A, \tilde{\rho})$ is a
left $H$-comodule
algebra over the $R$-bialgebroid $H$.

If $f:\ (A,{}^A\rho)\to (B,{}^B\rho)$ is a morphism of left $H$-comodule
algebras, then $f:\ (A,s_A,\widetilde{\rho})\to (B,s_B,\widetilde{\rho})$
is also a morphism of left $H$-comodule algebras over the $R$-bialgebroid
$H$. In view of the definition of $\widetilde{\rho}$, the left $H$-colinearity
is obvious. We also know that $f$ is an algebra map, so $f$ is a map of
$R$-rings
if $s_B=f\circ s_A$. Using the $H$-colinearity of $f$ and the fact that
$f(1_A)=1_B$, we find
\begin{eqnarray*}
f(s_A(\pi^L(h))&=&
\eps(1_{A<-1>}h)f(1_{A<0>})=\eps(f(1_A)_{<-1>}h)f(1_A)_{<0>}\\
&=& \eps(1_{B<-1>}h)f(1_{B<0>})=s_B(\pi^L(h))
\end{eqnarray*}

2) Conversely,
let $(A, s_A, {}^A\rho)$ be a left comodule algebra over the bialgebroid
$H$ as in
Definition~\ref{lhca}. We prove that $A$ is a weak left
$H$-comodule algebra with coaction given by
$\tilde{\rho} = \sigma \circ {}^A\rho$. Explicitly
$$
\tilde{\rho} (a) = a_{<-1>} \cdot S(1_{(1)}) \otimes 1_{(2)} \cdot a_{<0>} =
1_{(1)} a_{<-1>} \otimes 1_{(2)} \cdot a_{<0>} =
1_{(1)} a_{<-1>} \otimes s_A (1_{(2)}) a_{<0>},
$$
where we used that $H$ is a right $R$-module via the target map
$t_H = S^{-1}\mid_R$.
 The fact that $(A, \tilde{\rho})$ is a left
comodule of a weak Hopf algebra $H$ can easily be established
with the help of equations (\ref{weak.unit}).
We prove now that $\tilde{\rho}$ is an algebra map.
First, since ${\rm Im} ({}^A\rho) \subseteq H\times_R A$ we have
for all $r\in R$ and $a\in A$,
$$
a_{<-1>} t_H(r) \otimes_R a_{<0>} = a_{<-1>} \otimes_R a_{<0>} s_A(r).
$$
Applying the section $\sigma$ we obtain
\begin{equation}\label{joi1}
1_{(1)} a_{<-1>} t_H (r) \otimes s_A(1_{(2)}) a_{<0>} =
1_{(1)} a_{<-1>} \otimes s_A(1_{(2)}) a_{<0>} s_A(r).
\end{equation}
On the other hand, application of  $\sigma$ to an
expression reflecting the fact that
${}^A\rho : A \to H\times_R A$
is an algebra map leads to equality
\begin{equation}\label{joi2}
1_{(1)} (ab)_{<-1>} \otimes s_A(1_{(2)}) (ab)_{<0>} =
1_{(1)} a_{<-1>} b_{<-1>} \otimes s_A (1_{(2)})  a_{<0>} b_{<0>}
\end{equation}
Noting that $t_H =S^{-1}\mid_R$ and writing $1\sw{1'}\otimes 1\sw{2'}$
for another copy of $1\sw 1\otimes 1\sw 2$ we have
\begin{eqnarray*}
\tilde{\rho} (a) \tilde{\rho} (b)
&=& 1_{(1)} a_{<-1>} 1_{(1')} b_{<-1>} \otimes
s_A(1_{(2)}) a_{<0>} s_A(1_{(2')}) b_{<0>}\\
~~\mbox{(by Eq.\ (\ref{joi1}))}~~
&=& 1_{(1)} a_{<-1>} S^{-1} (1_{(2')}) 1_{(1')}b_{<-1>} \otimes
s_A(1_{(2)}) a_{<0>}b_{<0>}\\
~~\mbox{(by Eq.\ (\ref{joi2}))}~~
&=& 1_{(1)} (ab)_{<-1>} \otimes s_A(1_{(2)})(ab)_{<0>}
= \tilde{\rho} (ab),
\end{eqnarray*}
i.e., $\tilde{\rho}$ is multiplicative as required. It remains to prove
the unit property of the weak coaction $\tilde{\rho}$. Since
$\tilde{\rho} (1) = 1_{(1)} 1_{<-1>} \otimes 1_{(2)} \cdot 1_{<0>}$
we obtain that
$$
(H\otimes \tilde{\rho})\circ \tilde{\rho} (1) =
1_{(1)} 1_{<-2>} \otimes 1_{(2)} 1_{<-1>}\otimes 1_{(3)} \cdot 1_{<0>},
$$
hence the unit property of $\tilde{\rho}$ is equivalent to the
 following equation
$$
1_{(1')} \otimes 1_{(1)} 1_{<-1>} 1_{(2')} \otimes 1_{(2)} \cdot 1_{<0>} =
1_{(1)} 1_{<-2>} \otimes 1_{(2)} 1_{<-1>}\otimes 1_{(3)} \cdot 1_{<0>}.
$$
It is known, however, that
\begin{equation}\label{Stef}
1_{<-1>} \otimes_R 1_{<0>} = 1_H \otimes_R 1_A
\end{equation}
for $(A,{}^A\rho)$ is a
comodule algebra of an $R$-bialgebroid. Application of
$f$  yields
\begin{equation}
    \label{Gigel}
1_{(1)} 1_{<-1>} \otimes 1_{(2)}\cdot 1_{<0>} =
1_{(1)} \otimes s_A (1_{(2)}).
\end{equation}
Next apply $\Delta \otimes A$ to the preceding  equality to obtain
$$
1_{(1)} 1_{<-2>} \otimes 1_{(2)} 1_{<-1>}  \otimes 1_{(3)}\cdot 1_{<0>} =
1_{(1)} \otimes 1_{(2)} \otimes s_A (1_{(3)}).
$$
Now the required condition  follows from these two equations,
and equation (\ref{weak.unit}).
Thus we conclude that $(A, \tilde{\rho})$ is a left
weak $H$-comodule algebra.

Suppose now that $f:\ (A,s_A,{}^A\rho)\to (B,s_B,{}^B\rho)$ is a morphism
of left $H$-comodule algebras over the $R$-bialgebroid $H$. Then $f$ is
an algebra map, and
$$a_{<-1>}\otimes_Rf(a_{<0>})=f(a)_{<-1>}\otimes_R f(a)_{<0>}.$$
Applying $\sigma$ to both sides, we obtain
$$a_{<-1>}\cdot S(1_{(1)})\otimes 1_{(2)}\cdot f(a_{<0>})=
f(a)_{<-1>}\cdot S(1_{(1)})\otimes 1_{(2)}\cdot f(a)_{<0>}.$$
$f$ is a map of left $R$-modules, so the left hand side equals
$a_{<-1>}\cdot S(1_{(1)})\otimes  f(1_{(2)}\cdot a_{<0>})$, and this
means that $f:\ (A,\widetilde{\rho})\to (B,\widetilde{\rho})$
is left $H$-colinear, and $f$ is a morphism of left $H$-comodule 
algebras.\\

3) We still need to show that the functors constructed in parts 1) and 2)
of the proof are inverses to each other. First, let
$(A,\rho)$ be a left weak $H$-comodule algebra. It is first transformed into
a left $H$-comodule algebra $(A,\widetilde{\rho},s_A)$ over the
bialgebroid $R$, and then into a left weak $H$-comodule algebra
$(A,\ol{\rho})$.
We easily compute that
\begin{eqnarray*}
\ol{\rho}(a)&=& a_{<-1>}\cdot S(1_{(1)})\otimes 1_{(2)}\cdot a_{<0>}
= 1_{(1)}a_{<-1>}\otimes s_A(1_{(2)})a_{<0>}\\
&=& 1_{(1)}a_{<-1>}\otimes \epsilon(1_{<-1>}1_{(2)})1_{<0>}a_{<0>}
= 1_{<-2>}a_{<-1>}\otimes \epsilon(1_{<-1>})1_{<0>}a_{<0>}\\
&=& 1_{<-1>}a_{<-1>}\otimes 1_{<0>}a_{<0>}
= a_{<-1>}\otimes a_{<0>}=\rho(a),
\end{eqnarray*}
as needed. We used Definition \ref{def.weakDK}~(1).

Conversely, we start with a left $H$-comodule algebra $(A,\rho,s_A)$
over the bialgebroid $R$,
transform it into a weak left $H$-comodule algebra $(A,\widetilde{\rho})$
using part 2), and then back into $(A,\ol{\rho},\ol{s}_A)$ over $R$.
We have to show that $\rho=\ol{\rho}$ and $s_A=\ol{s}_A$.
We write $\rho(a)=a_{<-1>}\otimes_R a_{<0>}$, and then easily find that
\begin{eqnarray*}
\ol{\rho}(a)& = & a_{<-1>}\cdot S(1_{(1)})\otimes_R 1_{(2)}\cdot 
a_{<0>}\\
&= &a_{<-1>}\cdot (S(1_{(1)}) 1_{(2)})\otimes_R a_{<0>}
= a_{<-1>}\otimes_R a_{<0>}=\rho(a).
\end{eqnarray*}
Using equation (\ref{Gigel}), we obtain
$$
\ol{s}_A(\Pi^L(h)) = \eps(1_{(1)}1_{<-1>}h)s_A(1_{(2)})1_{<0>}=
s_A(\eps(1_{(1)}h)1_{(2)})=s_A(\Pi^L(h)).
$$
\end{proof}

\begin{proposition}\label{modulecoalgebra}
Let $H$ be a weak Hopf algebra viewed as an $R$-bialgebroid
as in Proposition \ref{webial}. The categories of left weak $H$-module
coalgebras (in the sense of
Definition~\ref{def.weakDK}) and of
left module coalgebras over the bialgebroid $H$ (in the sense of
Definition~\ref{altadef}) are isomorphic to each other.
\end{proposition}

\begin{proof}
1) Let $(C, \Delta_C, \eps_C)$ be a left weak $H$-module coalgebra and define
$\tilde{\Delta}_C : C\to C\otimes_R C$, $
c\mapsto  c_{(1)} \otimes_R c_{(2)}$ and  $\tilde{\eps}_C: C\to R$,
$c\mapsto \eps_C ( 1_{(1)} \cdot c) 1_{(2)}$. We will show that
these maps make $C$ into a left module coalgebra of the $R$-bialgebroid $H$.

Since $\Delta_C$
is left $H$-linear, so is
$\tilde{\Delta}_C$. This implies that $\tilde{\Delta}_C$ is an
$R$-bimodule map.
Next we prove that  $\tilde{\eps}_C$
is left $H$-linear. First note that setting $g=1_H$ in
Definition~\ref{def.weakDK}~(2) one immediately obtains $\eps_C(h\cdot c)
= \eps_H(h1\sw 2)\eps_C(1\sw 1\cdot c)$ for all $h\in H$, $c\in C$. In
particular
\begin{eqnarray*}
\tilde{\eps}_C (h\cdot c)
&=& \eps_C (1_{(1')} \cdot c) \eps_H ( 1_{(1)}h 1_{(2')} )1_{(2)}
= \Pi^L (h 1_{(2)}) \eps_C (1_{(1)} \cdot c) \\
&=& \tilde{\eps}_H (h 1_{(2)}) \eps_C (1_{(1)} \cdot c)
= \tilde{\eps}_H \Bigl(h \Pi^L ( \tilde{\eps}_C (c)  )   \Bigl)
= h \triangleright \tilde{\eps}_C (c),
\end{eqnarray*}
where $1\sw{1'}\otimes 1\sw{2'}$ is another copy of $\Delta(1_H)$ and
$\tilde{\eps}_H: H\to R$ is the counit of $H$ as an $R$-bialgebroid.
This implies that $\tilde{\eps}_C$ is an $R$-bimodule map.
Clearly, $\tilde{\Delta}_C$ is a coproduct and $\tilde{\eps}_C$ is a
counit of an
$R$-coring $C$.
The compatibility of $\tilde{\Delta}_C$ with the left action of $H$ on
$C$ follows immediately from the fact that $C$ is a weak module
coalgebra. Thus we conclude that $C$ is a left module coalgebra over the
$R$-bialgebroid $H$ as claimed.

Let $f:\ (C,\Delta_C,\eps_C)\to (D,\Delta_D,\eps_D)$ be a morphism of
left weak $H$-module coalgebras. Then $f$ is left $H$-linear, and it
clearly preserves the comultiplication over $R$. Furthermore
$$
\widetilde{\eps}_D(f(c))=\eps_D(1_{(1)}\cdot f(c))1_{(2)}=
\eps_D(f(1_{(1)}\cdot c))1_{(2)}
= \eps_C(1_{(1)}\cdot c)1_{(2)}=\widetilde{\eps}_C(c)
$$
and we conclude that $f:\ (C,\widetilde{\Delta}_C,\widetilde{\eps}_C)\to
(D,\widetilde{\Delta}_D,\widetilde{\eps}_D)$ is a morphism of left $H$-module
coalgebras over the $R$-bialgebroid $H$.\\

2) Conversely, assume that $C$ is a left module coalgebra over the
$R$-bialgebroid $H$,
in the sense of Definition~\ref{altadef}. We claim that $C$ is a left weak
$H$-module coalgebra with the coproduct
$$
\tilde{\Delta}_C : C\to C\otimes C, \quad
\tilde{\Delta}_C(c) = c_{(1)} \cdot S(1_{(1)}) \otimes 1_{(2)}\cdot c_{(2)} =
1_{(1)} \cdot c_{(1)} \otimes 1_{(2)}\cdot c_{(2)}
$$
(note that here $1\sw 1\otimes 1\sw 2 = \Delta(1_H)$ while $ c\sw
1\otimes _R c\sw 2 = \Delta_C(c)$), and the counit
$\tilde{\eps}_C :C\to k$, $c\mapsto \varphi (\eps_C (c) )$,
where $\eps_C : C\to R$ is the counit of the $R$-coring $C$ and
$ \varphi: R \to k$ is the Frobenius map $\varphi := \eps{|_{R}}$,
the restriction of the counit of a weak Hopf algebra $H$ to
the Frobenius algebra $R$.
This claim can be proven by a fairly straightforward calculation and hence
details of the proof are left to the reader.

Let $f:\ (C,\Delta_C,\eps_C)\to (D,\Delta_D,\eps_D)$ be a morphism of
left $H$-module coalgebras over the bialgebroid $H$. Then $f$ is left $H$-linear 
and
$$
(f\ot f)\widetilde{\Delta}_C(c)=
f(c_{
(1)}\cdot S(1_{(1)}))\ot f(1_{(2)}\cdot c_{(2)})
= f(c_{(1)})\cdot S(1_{(1)})\ot 1_{(2)}\cdot f(c_{(2)})
= \widetilde{\Delta}_D(f(c))
$$
and
$$\widetilde{\eps}_D(f(c))=\varphi(\eps_D(f(c)))=\varphi(\eps_C(c))=
\widetilde{\eps}_C(c),$$
so that $f$ is also a morphism in the category of weak left $H$-module
coalgebras.\\

3) We finally prove that the functors constructed in parts 1 and 2 of the
proof are inverses to each other. First we take a left weak $H$-module coalgebra
$(C,\Delta_C,\eps_C)$, turn it into a left module coalgebra
$(C,\widetilde{\Delta}_C,\widetilde{\eps}_C)$ over the bialgebroid $H$,
and then back into a left weak $H$-module coalgebra
$(C,\ol{\Delta}_C,\ol{\eps}_C)$. Using Remark \ref{separableremark},
we find
$$
\ol{\Delta}_C(c)= c_{(1)}\cdot S(1_{(1)})\ot 1_{(2)}\cdot c_{(2)}
= (\sigma\circ {\rm can})(c_{(1)}\ot c_{(2)})={\Delta}_C(c)
$$
and
$$\ol{\eps}_C(c)=\eps(\widetilde{\eps}_C(c))=\eps_C(1_{(1)}\cdot c)
\eps(1_{(2)})=\eps_C(c).$$
Finally take a left module coalgebra $(C,\Delta_C,\eps_C)$ over
the bialgebroid $H$, make it into a weak left $H$-module coalgebra
$(C,\widetilde{\Delta}_C,\widetilde{\eps}_C)$, and then back into
a left module coalgebra $(C,\ol{\Delta}_C,\ol{\eps}_C)$
over the bialgebroid $H$. Obviously
$$\ol{\Delta}_C(c)=c_{(1)}\cdot S(1_{(1)})\ot_R 1_{(2)}\cdot c_{(2)}=
\Delta_C(c).$$
Proving $\ol{\eps}=\eps$ is slightly more complicated. Recall that
$\eps:\ H\to k$ is the counit of the weak Hopf algebra $H$, and
$\eps_H=\Pi^L:\ H\to R$ is the counit of the $R$-bialgebroid $H$.
Take $c\in C$, and write $\eps_C(c)=\Pi^L(g)\in R$, for some $g\in H$.
Observe that $\eps(\Pi^L(g))=\eps(g)$ and 
$\widetilde{\eps}_C(c)=\eps(\eps_C(c))$, and compute
\begin{eqnarray*}
\ol{\eps}_C(c)&=& \widetilde{\eps}_C(1_{(1)}\cdot c)1_{(2)}
= \eps(\eps_C(1_{(1)}\cdot c))1_{(2)}\\
{\rm (by~Eq.~(\ref{Stef3}))}~~~&=& \eps(1_{(1)}\triangleright \eps_C(c))1_{(2)}\\
{\rm (by~Eq.~(\ref{Stef2}))}~~~&=& \eps(\eps_H(1_{(1)}s_H(\eps_C(c))))1_{(2)}\\
&=& \eps(\Pi^L(1_{(1)}\eps_C(c)))1_{(2)}
= \eps(1_{(1)}\eps_C(c))1_{(2)}\\
&=& \eps(1_{(1)}\Pi^L(g))1_{(2)}
= \eps(1_{(1)}\eps(1'_{(1)}g)1'_{(2)})1_{(2)}\\
{\rm (by~Eq.~(\ref{weak.unit}))}~~~&=&\eps(\eps(1_{(1)}g)1_{(2)})1_{(3)}\\
&=& \eps(1_{(1)}g)1_{(2)}=\Pi^L(g)=\eps_C(c)
\end{eqnarray*}
this completes the proof.
\end{proof}

Combining Propositions \ref{comodulealgebra} and \ref{modulecoalgebra},
we immediately obtain the following theorem, which is the main result
of this Section. We leave it to the reader to define morphisms
between left-left weak Doi-Koppinen data, and between left-left
Doi-Koppinen data over an algebra $R$.

\begin{theorem}\label{thm.weak-r}
Let $H$ be a weak Hopf algebra, and view it also as an $R$-bialgebroid,
as in Proposition \ref{webial}. Then there  is an isomorphism of
categories $\mathcal W\mathcal D\mathcal K(H)\cong \mathcal D\mathcal K_R(H)$.
Furthermore, the corresponding categories of Doi-Koppinen modules are
isomorphic.
\end{theorem}

\section{A coring associated to a Doi-Koppinen datum over $R$ and
 applications}
In this section we construct an $A$-coring corresponding to a given
Doi-Koppinen datum $(H,A,C)_R$ over $R$. This allows
one to use methods employed in \cite{Brz:cor} to derive various
properties of Doi-Koppinen modules over an algebra.

\begin{proposition}\label{prop.coring}
Let $(H,A,C)_R$ be a Doi-Koppinen datum over $R$. Then
$\cC = C\otimes_R A$ is an $A$-bimodule with the right action given by
the multiplication in $A$ and the left action $a\cdot(c\otimes_Ra') =
a\sco{-1}\cdot c\otimes_R a\sco 0a'$, for all $a,a'\in A$, $c\in C$.
Furthermore $\cC$ is an $A$-coring
with comultiplication $\Delta_\cC = \Delta_C\otimes_R A$ and the counit
$\eps_\cC = \eps_C\otimes_R A$, where $\Delta_C$, $\eps_C$ are the
coproduct and the counit of the $R$-coring $C$. In this case the
categories of left $\cC$-comodules and of left-left Doi-Koppinen
modules over $R$ are isomorphic to each other.
\end{proposition}
\begin{proof}
First note that the left action of $A$ on $\cC$ is well-defined since
the image of the left $H$-coaction of $A$ is required to be in
$A\times_R H$. The fact that it is an action indeed follows from the
fact that $A$ is a left $H$-comodule algebra.
    Note also that in the definitions of $\Delta_\cC$ and
$\eps_\cC$ we used the natural isomorphisms $C\otimes_R A\otimes_A
C\otimes _R A\cong C\otimes_R C\otimes_R A$ and $R\otimes_R A\cong A$
respectively. Clearly $\Delta_\cC$ is a right $A$-module map. To prove
that it is a left $A$-module map as well take any $a,a'\in A$ and $c\in
C$ and compute
\begin{eqnarray*}
\Delta_\cC(a\cdot (c\otimes_R a')) &=&  (a\sco{-1}\cdot c)\sw 1\otimes_R
(a\sco{-1}\cdot c)\sw 2\otimes_R a\sco 0a'\\
& =& a\sco{-2}\cdot c\sw 1\otimes_R
a\sco{-1}\cdot c\sw 2\otimes_R a\sco 0a',
\end{eqnarray*}
where we used that $C$ is a left $H$-module coalgebra. On the other
hand
\begin{eqnarray*}
a\cdot\Delta_\cC(c\otimes_R a') &=& a\cdot(c\sw 1\otimes_R
1)\otimes_A(c\sw 2\otimes_R a') = a\sco{-1}\cdot c\sw 1\otimes_R a\sco
0\cdot (c\sw 2\otimes_R a')\\
&=& a\sco{-2}\cdot c\sw 1\otimes_R
a\sco{-1}\cdot c\sw 2\otimes_R a\sco 0a'.
\end{eqnarray*}
This proves that $\Delta_\cC$ is right $A$-linear, hence it is an
$A$-bimodule map as required. Directly from the definition of
$\Delta_{\cC}$ it follows that it is coassociative.

It is clear that $\eps_\cC$ is right $A$-linear. To prove that it is also
a left $A$-module morphism take any $a,a'\in A$, $c\in C$ and compute
\begin{eqnarray*}
\eps_\cC(a\cdot(c\otimes_R a')) & = &\eps_C(a\sco{-1}\cdot c)\cdot (a\sco
0a')
= \eps_H(a\sco{-1}s_H(\eps_C(c)))\cdot (a\sco 0a')\\
& = & \eps_H(a\sco{-1}t_H(\eps_C(c)))\cdot (a\sco
0a')= \eps_{H}(a\sco{-1})\cdot (a\sco 0 s_{A}(\eps_{C}(c))a')\\
&=& a\sco 0 s_{A}(\eps_{C}(c))a' = a\eps_{\cC}(c\otimes_{R}a'),
\end{eqnarray*}
where we used that $C$ is a left $H$-module coalgebra to obtain the
second equality, then equation (\ref{25}) to derive the third one
and the fact that the image
of the left coaction of $H$ on $A$ is in $H\times_{R}A$ to obtain the
fourth equality. This proves that $\eps_{\cC}$ is left $A$-linear
hence it is $A$-bilinear. The fact that $\eps_{\cC}$ is a counit of
$\cC$ follows directly from the definition of $\eps_{\cC}$. Thus we
conclude that $\cC$ is an $A$-coring as stated.

To prove the isomorphism of categories, take any left $\cC$-comodule
$(M, \rho)$ and view it as a Doi-Koppinen module via the
same coaction $\rho: M\to C\otimes_{R}A\otimes_{A}M \cong
C\otimes_{R}M$. Conversely, any Doi-Koppinen module $(M,\cdot,
\rho)$ can be viewed as a left $\cC$ comodule via $\rho:
M\to C\otimes_{R}M\cong C\otimes_{R}A\otimes_{A}M = \cC\otimes_{A}M$.
\end{proof}

One can now use the general results about corings in
\cite{Brz:cor}\footnote{Note, however, that some care
has to be taken when applying \cite{Brz:cor} since
this paper is formulated in the right-right module convention.}
combined with
Proposition~\ref{prop.coring} to derive various properties of
Doi-Koppinen modules over a noncommutative algebra $R$. For example
\cite[Lemma~3.1]{Brz:cor} implies that the forgetful functor $F:
{}^C_A\M (H;R) \to {}_{A}\M$ is the left adjoint of the induction
functor $G = C\otimes_{R}-: {}_{A}\M\to {}^{C}_{A}\M(H;R)$. Furthermore by
\cite[Theorems~3.3,~3.5]{Brz:cor} one has
\begin{corollary}
Let $(H,A,C)_R$ be a left-left Doi-Koppinen datum over $R$.

(1) The induction
functor $G= C\otimes_{R}-: {}_{A}\M\to {}^{C}_{A}\M(H;R)$ is separable if
and only
if there exists $e=\sum_{i}c^{i}\otimes_{R}a^{i}\in C\otimes_{R}A$
such that $\sum_{i}\eps_{C}(c^{i})\cdot a^{i}=1_{A}$ and for all $a\in
A$, $\sum_{i}a\sco{-1}\cdot c^{i}\otimes_{R}a\sco 0a^{i}=
\sum_{i}c^{i}\otimes_{R}a^{i}a$.

(2) The forgetful functor $F:{}^C_A\M (H;R) \to {}_{A}\M$ is
separable if and only if there exists a right $R$-bimodule map
$\gamma: C\otimes_{R}C\to A$ such that for all $a\in A$ and $c,c'\in C$,
\begin{itemize}
\item $\gamma(c\sw 1\otimes_{R}c\sw 2) = \eps_{C}(c)\cdot 1_{A}$,
\item $\gamma(a\sco{-2}\cdot c\otimes_{R}a\sco{-1}\cdot c')a\sco 0 =
a\gamma(c\otimes_{R}c')$,
\item
 $c\sw 1\otimes_{R}\gamma(c\sw
2\otimes_{R}c') = \gamma(c\otimes_{R}c'\sw 1)\sco{-1}\cdot c'\sw
2\otimes _{R}\gamma(c\otimes_{R}c'\sw 1)\sco 0$.
\end{itemize}
\end{corollary}

Finally, it has been observed in \cite[Proposition~2.3]{Brz:cor},
that given a weak entwining structure $(A,C,\psi)$, and hence a weak
 Doi-Koppinen
datum in particular, one can construct a coring obtained as an image
of certain projection in $C\otimes A$. Since a weak Doi-Koppinen datum
is a special case of a Doi-Koppinen datum over $R$ it is important to
study the relationship of this coring to $C\otimes_{R}A$.

\begin{proposition}
    Let $(H,A,C)$ be weak Doi-Koppinen datum. Define the
    corresponding $A$-coring $\tilde{\cC} =
    \{\sum_{i}1\sco{-1}\cdot c^{i}\otimes 1\sco 0a^{i}\; | \; a^{i}\in
    A, c^{i}\in C\}$, with the coproduct $\Delta_{\tilde{\cC}} =
    (\Delta_{C}\otimes A)\mid_{\tilde{\cC}}$ and the counit
    $\eps_{\tilde{\cC}} = (\eps_{C}\otimes A)\mid_{\tilde{\cC}}$.
    View $(H,A,C)$ as a Doi-Koppinen datum over $R={\rm Im}\Pi^L$ in
    Theorem~\ref{thm.weak-r}. Then $\tilde{\cC} \cong \cC =
    C\otimes_{R}A$ as $A$-corings.
\end{proposition}
\begin{proof}
    Consider two maps $\theta:\cC\to \tilde{\cC}$,
    $c\otimes_{R}a\mapsto 1\sco{-1}\cdot c\otimes 1\sco 0a$, and
    $\tilde{\theta}: \tilde{\cC}\to \cC$, $\sum_{i}(c^{i}\otimes
    a^{i})\mapsto \sum_{i}c^{i}\otimes_{R}a^{i}$. The map $\theta$ is
    well-defined because using equation~(\ref{eq.sep}) we have for all
    $a\in A$, $c\in C$ and $r\in R ={\rm Im}\Pi^{L}$,
    \begin{eqnarray*}
	\theta(c\cdot r\otimes _{R}a) &= &1\sco{-1}\cdot (c\cdot r)\otimes
	1\sco 0a=1\sco{-1}t_{H}(r)\cdot c\otimes 1\sco 0a\\
	&=& 1\sco{-1}\cdot c\otimes 1\sco 0s_{A}(r)a =
	\theta(c\otimes_{R}r\cdot a).
    \end{eqnarray*}
    Clearly, $\theta$ is a right $A$-module map. Recall from
    \cite[Proposition~2.3]{Brz:cor} that $\tilde{\cC}$ is a left
    $A$-module via $a\cdot\left(\sum_{i}1\sco{-1}\cdot c^{i}\otimes 1\sco
    0a^{i}\right) = \sum_{i}a\sco{-1}\cdot c^{i}\otimes a\sco 0 a^{i}$. Now,
    the fact that $A$ is a weak left $H$-comodule algebra implies
    that $\theta$ is a left $A$-module map. Thus $\theta$ is an
    $A$-bimodule map. To prove that $\theta$ is a coring map take any
    $a\in A$, $c\in C$ and use the fact that $A$ is an $H$-comodule
    algebra to compute
    \begin{eqnarray*}
	(\theta\otimes_{A}\theta)\circ\Delta_{\cC}(c\otimes_{R}a) & = &
    1\sco{-2}\cdot c\sw 1\otimes(1\sco{-1}1\sco{-1'})\cdot c\sw
    2\otimes1\sco 01\sco{0'}a\\
    &= &1\sco{-2}\cdot c\sw 1\otimes1\sco{-1}\cdot c\sw
    2\otimes1\sco 0a,
    \end{eqnarray*}
    where $1\sco{-1'}\otimes1\sco{0'}$  is another copy of
    $1\sco{-1}\otimes 1\sco 0$. On the other hand, using the fact
    that $C$ is a left $H$-module coalgebra we have
    $$
    \Delta_{\tilde{\cC}}(\theta(c\otimes_{R}a)) = (1\sco{-1}\cdot
    c)\sw 1\otimes (1\sco{-1}\cdot c)\sw 2 \otimes 1\sco 0a =
    1\sco{-2}\cdot c\sw 1\otimes1\sco{-1}\cdot c\sw
    2\otimes1\sco 0a,
    $$
    as required. Thus we have proven that $\theta$ is a map of
    $A$-corings. We now prove that $\tilde{\theta}$ is an inverse of
    $\theta$. For a typical element
    $x=\sum_{i}1\sco{-1}\cdot c^{i}\otimes 1\sco 0a^{i}$ of $\tilde{\cC}$
    we have
    $$\theta\circ\tilde{\theta}(x) = \sum_{i}1\sco{-1'}1\sco{-1}\cdot
    c^{i}\otimes 1\sco{0'}1\sco 0a^{i} = \sum_{i}1\sco{-1}\cdot
    c^{i}\otimes 1\sco 0a^{i}= x,
    $$
    for $A$ is a weak $H$-module algebra. On the other hand, since
    $1\sco{-1}\otimes 1\sco 0 = 1\sw 1\otimes s_A(1\sw 2)$ (cf.\ proof
    of Theorem~\ref{thm.weak-r}) we have for all $a\in A$, $c\in C$
    \begin{eqnarray*}
    \tilde{\theta}\circ\theta(c\otimes_{R}a) & = &
	1\sco{-1}\cdot c\otimes
    _{R}1\sco 0a = 1\sw 1\cdot c\otimes_{R}s_{A}(1\sw 2)a \\
    & = &(1\sw
    1\cdot c)\cdot 1\sw 2\otimes_{R}a = S^{-1}(1\sw 2) 1\sw 1\cdot
    c\otimes _{R}a = c\otimes_{R}a.
    \end{eqnarray*}
    This completes the proof that $\theta$ is an isomorphism of
    $A$-corings.
\end{proof}

\begin{center}
{\sc Acknowledgements}
\end{center}
T.\ Brzezi\'nski thanks
EPSRC for an Advanced Research Fellowship, S.\ Caenepeel thanks
Department of Mathematics, University of Wales Swansea for
hospitality, and G.\ Militaru thanks the Royal Society for a visiting
fellowship.

\end{document}